\documentclass[12pt]{article}
\usepackage{amsfonts}
\usepackage{}
\usepackage{amsthm,amsmath}
\usepackage{amscd}
 \usepackage{amssymb}
 \usepackage{hyperref}
 \usepackage[all]{xypic}
 \usepackage{mathtools}
  \usepackage[numbers,sort&compress]{natbib}
 \usepackage[utf8]{inputenc}
  \newtheorem{lemma}{Lemma}[section]
 \newtheorem{corollary}[lemma]{Corollary}
 
 \newtheorem{theorem}[lemma]{Theorem}

\usepackage[utf8]{inputenc} 

\def\includegraphics{}

\usepackage{amscd}
\usepackage{amsmath,amssymb,mathrsfs}
\usepackage{hyperref}
\usepackage[all]{xypic}
\usepackage{times}
\usepackage{inputenc}
\usepackage{titlesec}
 \usepackage{amsthm}

\advance\headheight1.15pt
\newenvironment{proof of Theorem 1.1 step1}{{\noindent\it {\bf Proof of the implication $(1)\Rightarrow(3)$}}\quad}{\hfill $\square$\par}
\newenvironment{proof of Theorem 1.1 step2}{{\noindent\it {\bf Proof of the implication $(3)\Rightarrow(2)$}}\quad}{\hfill $\square$\par}
\newenvironment{proof of Theorem 1.2}{{\noindent\it {\bf Proof of Theorem 1.2}}\quad}{\hfill $\square$\par}
\newcommand{\I}{\int_{0}^{1}}      
   
\newcommand{\sn}{\sum_{n=0}^{\infty}}           
\newcommand{\hd}{H(\mathbb{D})}            
\newcommand{\dd}{\mathbb{D}}           

\newcommand{\Cu}{\mathcal{C}_\mu}   \newcommand{\cua}{\mathcal{C}_{\mu,\alpha}}

\newcommand{\comment}[1]{}

\usepackage{lineno}

\begin{document}

\baselineskip=8pt
\title{\fontsize{15}{0}\selectfont  Generalized Ces\`{a}ro-like operator from a class of analytic function spaces to analytic Besov spaces}
\author{\fontsize{11}{0}\selectfont
   Pengcheng Tang$^{*,a}$ \\
   \fontsize{10}{0}\it{
 $^{a}$ Hunan University of Science and Technology, Xiangtan, Hunan 411201, China} 
}
\date{}
\maketitle
\thispagestyle{empty}
\begin{center}
\textbf{\underline{ABSTRACT}}
\end{center}
\ \ \ \ Let $\mu$ be a   finite  positive Borel measure  on $[0,1)$ and $f(z)=\sum_{n=0}^{\infty}a_{n}z^{n} \in H(\mathbb{D})$. For $0<\alpha<\infty$, the generalized Ces\`aro-like operator $\cua$ is defined by
$$
\mathcal {C}_{\mu,\alpha}(f)(z)=\sum^\infty_{n=0}\left(\mu_n\sum^n_{k=0}\frac{\Gamma(n-k+\alpha)}{\Gamma(\alpha)(n-k)!}a_k\right)z^n, \ z\in \dd,
$$
where, for $n\geq 0$, $\mu_n$ denotes the $n$-th  moment of the measure  $\mu$, that is,
$\mu_n=\int_{0}^{1} t^{n}d\mu(t)$.

For $s>1$,  let $X$ be a Banach subspace of  $\hd$ with $\Lambda^{s}_{\frac{1}{s}}\subset X\subset\mathcal {B}$. In this paper,  for $1\leq p <\infty$, we  characterize the measure $\mu$ for which  $\cua$ is bounded(or compact) from $X$ into analytic Besov
space $B_{p}$.

\begin{flushleft}
\end{flushleft}
\begin{flushleft}
\end{flushleft}

\let\thefootnote\relax\footnote{$^*$Corresponding Author}
\let\thefootnote\relax\footnote{ Pengcheng Tang: www.tang-tpc.com@foxmail.com}

\vspace{1cm}
%

\section{Introduction} \label{Sec:Intro}

 \ \ \  \ \ \ Let $\mathbb{D}=\{z\in \mathbb{C}:\vert z\vert <1\}$ denote the open unit disk of the complex plane $\mathbb{C}$ and $H(\mathbb{D})$ denote the space of all
analytic functions in $\mathbb{D}$. $H^{\infty}$ denote the set  of bounded analytical functions on $\dd$.


The Bloch space $\mathcal {B}$ consists of those functions $f\in  H(\mathbb{D}) $ for which
$$
\vert \vert f\vert \vert _{\mathcal {B}}=\vert f(0)\vert +\sup_{z\in \mathbb{D}}(1-\vert z\vert ^{2})\vert f'(z)\vert <\infty.
$$

For $1<p<\infty$,   the analytic Besov space $B_p$  consists of those  functions $f\in \hd $ such that
$$
\|f\|_{B_p}=|f(0)|+\left(\int_{\dd}|f'(z)|^p(1-|z|)^{p-2}\mathrm{d}A(z)\right)^{\frac{1}{p}}<\infty,
$$
where $dA(z)=\frac{dxdy}{\pi}$ is the  normalized area measure on $\dd$.
When $p=2$, then $B_2$ is just the classic Dirichlet space.
If  $1<p_1<p_2<\infty$, then  $B_{p_1}\subsetneq B_{p_2} \subsetneq\mathcal {B}$. It is known that  the   analytic Besov spaces are M\"{o}bius invariant and the Bloch space   $\mathcal {B}$ is largest M\"{o}bius invariant  space.

The space $B_{1}$ consists of $f\in \hd$ such that
$$\int_{\dd}|f''(z)|dA(z)<\infty.$$
 We know that $B_{1}$ is the smallest M\"{o}bius invariant Banach spaces of analytic function in $\dd$ and $B_{1}\subsetneq H^{\infty}$. See \cite[Chapter 5]{H14} for the theory of these spaces.

Let  $1\leq p<\infty$ and $0<\alpha\leq 1$, the mean Lipschitz space $\Lambda^p_\alpha$ consists of those functions $f\in H(\dd)$ having a non-tangential limit  almost everywhere such that $\omega_p(t, f)=O(t^\alpha)$ as $t\to 0$. Here $\omega_p(\cdot, f)$ is the integral modulus of continuity of order $p$ of the function $f(e^{i\theta})$. It is  known (see \cite{b1}) that $\Lambda^p_\alpha$ is a subset of $H^p$ and
$$\Lambda^p_\alpha=\left(f\in H(\dd):M_p(r, f')=O\left(\frac{1}{(1-r)^{1-\alpha}}\right), \ \ \mbox{as}\ r\rightarrow 1\right).$$
The space $\Lambda^p_\alpha$ is a Banach space with the norm $||\cdot||_{\Lambda^p_\alpha}$ given by
$$
\|f\|_{\Lambda^p_\alpha}=|f(0)|+\sup_{0\leq r<1}(1-r)^{1-\alpha}M_p(r, f').
$$
It is known (see e.g. \cite{bu}) that
$$\Lambda^{p}_{\frac{1}{p}}\subsetneq \mathcal {B}. \ \ 1<p<\infty.$$

For $f(z)=\sum_{n=0}^\infty a_nz^n\in \hd$, the Ces\`{a}ro operator $\mathcal {C}$ is defined  by
 $$
\mathcal {C}(f)(z)=\sum_{n=0}^\infty\left(\frac{1}{n+1}\sum_{k=0}^n a_k\right)z^n, \  z\in\dd.
$$

\comment{
The Ces\`{a}ro operator   $\mathcal {C}$ is  bounded on $H^{p}$ for $0<p<\infty$. The case of $1<p<\infty$
follows from a result of Hardy on Fourier series \cite{ha} together with the Riesz
transform. Siskakis \cite{sis1} give an alternative proof of this result and to extend it to
$p = 1$ by using semigroups of  composition operators.  A direct proof of the boundedness on $H^1$ was given by Siskakis in \cite{sis2}. Miao
\cite{miao} proved the case $0 <p < 1$. Stempak \cite{ste} gave a proof valid for $ 0 < p \leq 2$.
 Andersen \cite{ad} and Nowak \cite{no} provided another proof valid for all $0<p<\infty$. In the case $p=\infty$, since $\mathcal {C}(1)(z)=\log\frac{1}{1-z}\notin H^{\infty}$, so that $\mathcal {C}(H^{\infty})\nsubseteq H^{\infty}$. Danikas  and  Siskakis \cite{sis4} proved  that $\mathcal {C}(H^{\infty})\nsubseteq BMOA$ and $\mathcal {C}(BMOA)\nsubseteq BMOA$. Ces\`{a}ro operator   $\mathcal {C}$ act on  weighted Bergman spaces, Dirichlet space  and general mixed normed spaces $H(p,q,\varphi)$ the reader is referred to \cite{sis2,v,ale,gla,shi}.
 }

 The boundedness and compactness of the Ces\`{a}ro operator $\mathcal {C}$ and its generalizations
on various spaces of analytic functions such as  Hardy spaces, Bergman spaces, Dirichlet spaces, Bloch space, $Q_{p}$ space, mixed norm space
have been widely studied. See  \cite{ces9,ces8,ces13,ces15,ces4,ces3,ces14,ces5}
 and the references therein.

Recently, Galanopoulos, Girela and Merch\'an \cite{ces1} introduced a Ces\`aro-like operator $\Cu$ on $\hd$, which is a natural generalization of the classical Ces\`{a}ro operator $\mathcal {C}$.
They consider the following generalization: For a positive Borel measure
$\mu$ on the interval $[0, 1)$ they define the operator
 $$
\Cu (f)(z)=\sum^\infty_{n=0}\left(\mu_n\sum^n_{k=0}\widehat{f}(k)\right)z^n=\int_{0}^{1}\frac{f(tz)}{(1-tz)}d\mu(t), \  z\in\dd. \eqno{(1.3)}
$$
where $\mu_{n}$ stands for the  moment of order $n$ of $\mu$, that is,
  $\mu_{n}=\I t^{n}d\mu(t)$. They studied the  operators $\Cu$  acting on distinct spaces of analytic functions(e.g. Hardy space, Bergman space, Bloch space, etc.).

The  Ces\`aro-like operator $\Cu$ defined above has  attracted the interest  of many mathematicians.
For instance,  Jin and Tang \cite{ces6} studied  the boundedness(compactness) of  $\Cu$  from one Dirichlet-type space $\mathcal {D}_{\alpha}$  into another one $\mathcal {D}_{\beta}$.
Bao, Sun and Wulan  \cite{baoo} studied the range of $\Cu$ acting on $H^{\infty}$.  Blasco \cite{blas} investigated the operators $\Cu$  induce by  complex Borel measures on $[0, 1)$,  and extended the results of  \cite{ces1}  to this more general case.
  The operators  $\Cu$ associated to arbitrary complex Borel measures on  $\dd$  the reader is referred to  \cite{04}.

Bao et al. \cite{baoo} introduced a more general  Ces\`{a}ro-like operator: Suppose that $0<\alpha<\infty$
and  $\mu$ is a finite positive Borel measure on $[0, 1)$.  For $f(z)=\sn a_{n}z^{n} \in \hd$, they defined
$$
\mathcal {C}_{\mu,\alpha}(f)(z)=\sum^\infty_{n=0}\left(\mu_n\sum^n_{k=0}\frac{\Gamma(n-k+\alpha)}{\Gamma(\alpha)(n-k)!}a_k\right)z^n, \ z\in \dd.
$$
A simple calculation with power
series gives the integral form of $\cua$ as follows.
$$\cua(f)(z)=\int_{0}^{1}\frac{f(tz)}{(1-tz)^{\alpha}}d\mu(t).$$
It is clear that $\mathcal {C}_{\mu,1}=\Cu$.

For $1<s<\infty$,  let  $X$ be a  Banach subspace of  $\hd$ with $\Lambda^{s}_{\frac{1}{s}}\subset X\subset\mathcal {B}$. There are many well known spaces located  between the mean Lipschitz space $\Lambda^{s}_{\frac{1}{s}}$ and
the Bloch space $\mathcal {B}$. 
 In \cite{baoo}, the authors investigated the range of $\cua$ acting on $H^{\infty}$. They proved that if $\max\{1,\frac{1}{\alpha}\}<s<\infty$, then $\cua(H^{\infty})\subset X$ if and only if $\mu$ is an $\alpha$-Carleson measure. Zhou \cite{05} considers the same problem for  the measure $\mu$ supported on $\dd$.  Galanopoulos  et al. \cite{03} studied  the behaviour of the operators $\mathcal {C}_{\mu,1}$  on the Dirichlet space and on the
analytic Besov spaces.
Sun et al. \cite{sun} studied the operator $\mathcal {C}_{\mu,1}$ acting from $B_{p}$ to $X$ recently.
 It remains open to characterize the boundedness and the compactness of $\cua$  from $B_{p}$ to $B_{p}$ $(p>1)$. The Besov spaces $B_{p}$ and Bloch space $\mathcal {B}$ are M\"{o}bius invariant and the Bloch space $\mathcal {B}$ can be regarded as the limit case of $B_{p}$ as $p\rightarrow +\infty$.
The purpose of  this paper is  describe the measure $\mu$ such that the operator $\cua$  is bounded(and compact) from $X$ to $B_{p}$, $1\leq p<\infty$.

 Our main results are included in the following.
 \begin{theorem}\label{th1.1}
Suppose $0<\alpha<\infty$, $1<s<\infty$,  $\max\{1,\frac{1}{\alpha}\}\leq p<\infty$. Let $\mu$ be a  finite positive  Borel measure on  $[0, 1)$ and $X$  is a Banach subspace of  $\hd$ with $\Lambda^{s}_{\frac{1}{s}}\subset X\subset\mathcal {B}$.  Then the following statements are equivalent.
\\ (1) The operator $\cua$ is bounded from $X$ to $B_{p}$.
\\(2)  The operator $\cua$ is compact  from $X$  to $B_{p}$.
\\(3) The measure $\mu$ satisfies
$$\sn (n+1)^{p\alpha-1}\mu_{n}^{p}\log^{p}(n+2)<\infty.$$
\end{theorem}

For $p = 1$, we have the following corollary.
\begin{corollary}
Suppose $1\leq \alpha<\infty$ and $1<s<\infty$. Let $\mu$ be a  finite positive  Borel measure on  $[0, 1)$ and $X$  is a Banach subspace of  $\hd$ with $\Lambda^{s}_{\frac{1}{s}}\subset X\subset\mathcal {B}$.  Then the following statements are equivalent.
\\ (1) The operator $\cua$ is bounded from $X$ to $B_{1}$.
\\(2)  The operator $\cua$ is compact  from $X$  to $B_{1}$.
\\(3) The measure $\mu$ satisfies
$$\sn (n+1)^{\alpha-1}\mu_{n}\log(n+2)<\infty.$$
\\(4) The measure $\mu$ satisfies
$$\int_{0}^{1}\frac{\log\frac{e}{1-t}}{(1-t)^{\alpha}}d\mu(t).$$
\end{corollary}

 Throughout the paper, the letter $C$ will denote an absolute constant whose value depends on the parameters
indicated in the parenthesis, and may change from one occurrence to another. We will use
the notation $``P\lesssim Q"$ if there exists a constant $C=C(\cdot) $ such that $`` P \leq CQ"$, and $`` P \gtrsim Q"$ is
understood in an analogous manner. In particular, if  $``P\lesssim Q"$  and $ ``P \gtrsim Q"$ , then we will write $``P\asymp Q"$.


\section{ Preliminaries} \label{prelim}

To prove our main results, we need some  preliminary results which will be repeatedly
used throughout the rest of the paper. We begin with a characterization of the functions  $f\in \hd$ whose sequence of Taylor coefficients is
decreasing which belong to  $B^{p}$.

\begin{lemma}\label{lem2.1}
Let $1<p<\infty$ and $f(z)=\sn a_{n}z^{n}\in \hd$.  Suppose that the sequence $\{a_{n}\}_{n=0}^{\infty}$ is a decreasing sequence of non-negative
real numbers. Then $f\in B_{p}$ if and only if
$$\sum_{n=1}^{\infty}n^{p-1}a_{n}^{p}<\infty.$$
\end{lemma}
This result can be proved with arguments similar to those used in the proofs of \cite[Theorem 3.1]{I1}.  For a detailed proof, see also \cite[Theorem 3.10]{H5}.

The following lemma  contains a characterization of $L^{p}$-integrability of power series
with nonnegative coefficients.  For a proof, see \cite[Theorem 1]{1983}.
\begin{lemma}\label{lem2.2}
Let  $0<\beta,p<\infty$, $\{\lambda_{n}\}_{n=0}^{\infty}$  be a sequence of  non-negative  numbers. Then
$$\int_{0}^{1}(1-r)^{p\beta-1}\left(\sn \lambda_{n}r^{n}\right)^{p}dr\asymp \sn 2^{-np\beta}\left(\sum_{k\in I_{n}}\lambda_{k}\right)^{p}$$
where  $I_{0}=\{0\}$, $I_{n}=[2^{n-1},2^{n})\cap \mathbb{N}$ for $n\in \mathbb{N}$.
\end{lemma}

The following lemma is a consequence of Theorem 2.31 on page 192 of the classical monograph \cite{b8}.
\begin{lemma}\label{lem2.3}
 (a)\  The Taylor coefficients $a_{n}$ of the function
 $$f(z)=\frac{1}{(1-z)^{\beta}}\log^{\gamma}\frac{2}{1-z},  \ \ \beta>0,\gamma\in \mathbb{R}, \ z\in \mathbb{D}$$
 have the property $a_{n}\asymp n^{\beta-1}(\log(n+1))^{\gamma}$.

 (b)\  The Taylor coefficients $a_{n}$ of the function
 $$f(z)=\log^{\gamma}\frac{2}{1-z}, \ \ \gamma >0, z\in \mathbb{D}$$
 have the property $a_{n}\asymp n^{-1}(\log(n+1))^{\gamma-1}$.
\end{lemma}

We  also need the following estimates (see, e.g. Proposition  1.4.10 in \cite{b3}).
\begin{lemma}\label{lem2.4}
 Let $\alpha$ be any real number and $z\in \dd$.  Then
$$
\int^{2\pi}_0\frac{d\theta}{|1-ze^{-i\theta}|^{\alpha}}\asymp
\begin{cases}1 & \enspace \text{if} \ \ \alpha<1,\\
                     \log\frac{2}{1-|z|^2} & \enspace  \text{if} \ \  \alpha=1,\\
                     \frac{1}{(1-|z|^2)^{\alpha-1}} & \enspace \text{if}\ \  \alpha>1,
                   \end{cases}
$$
\end{lemma}

The following lemma is  useful in dealing
with the compactness. The proof is similar to that of Proposition 3.11 in \cite{hb1}. The details are omitted.
\begin{lemma}\label{lem2.5}
Let $p\geq 1$, $s>1$, $X$  be a Banach subspace of  $\hd$ with $\Lambda^{s}_{\frac{1}{s}}\subset X\subset\mathcal {B}$. Suppose that $T$ is a bounded operator  from $X$ to $B_{p}$. Then $T$ is compact if and only if for any bounded sequence $\{f_{k}\}$ in $X$
which converges to $0$ uniformly on every compact subset of  $\dd$, we have $\lim_{k\rightarrow\infty}||T(f_{k})||_{B_{p}}=0$.
\end{lemma}

 \section{ Proofs of the main results} \label{prelim}

We now present the proofs of Theorem 1.1.

\begin{proof of Theorem 1.1 step1}
Since the definition of the space $B_{1}$ is slightly different from $B_{p}(p>1)$, we split the proof into $p=1$  and  $p>1$.

{\bf Case $p=1$.} \ \ Assume $\cua$ is bounded from $X$ to $B_{1}$. Let $g(z)=\log\frac{1}{1-z}=\sum^{\infty}_{k=1}\frac{z^{k}}{k}$,  it is easy to check that  $g\in \Lambda^{s}_{\frac{1}{s}}\subset X$. This implies that  $\cua(g)\in B_{1}$. For $z\in \dd$, by the definition of $\cua$ we get
$$\cua(g)''(z)=\sum_{n=0}^{\infty}\left((n+2)(n+1)\mu_{n+2}\sum_{k=1}^{n+2}\frac{\Gamma(n+2-k+\alpha)}{\Gamma(\alpha)(n+2-k)!k}\right)z^{n}.$$
For $0<r<1$, the Hardy's inequality shows that
$$M_{1}(r,\cua(g)'')\gtrsim \sn \left((n+2)\mu_{n+2}\sum_{k=1}^{n+2}\frac{\Gamma(n+2-k+\alpha)}{\Gamma(\alpha)(n+2-k)!k}\right)r^{n}.$$
Hence,
 \[ \begin{split}
 1&\gtrsim ||g||_{X}\gtrsim ||\cua(g)||_{B_{1}}=\int_{\dd}|\cua(g)''(z)|dA(z)\\
 & = 2\int_{0}^{1}M_{1}(r,\cua(g)'')rdr\\
 & \gtrsim \int_{0}^{1} \sn \left((n+2)\mu_{n+2}\sum_{k=1}^{n+2}\frac{\Gamma(n+2-k+\alpha)}{\Gamma(\alpha)(n+2-k)!k}\right)r^{n+1}dr\\
 & \gtrsim  \sn \mu_{n+2}\sum_{k=1}^{n+2}\frac{\Gamma(n+2-k+\alpha)}{\Gamma(\alpha)(n+2-k)!k}.
   \end{split} \]
 Using the Stirling's formula  we get
   $$\sum_{k=1}^{n+2}\frac{\Gamma(n+2-k+\alpha)}{\Gamma(\alpha)(n+2-k)!k}\asymp \sum_{k=1}^{n+2}\frac{(n+3-k)^{\alpha-1}}{k}. $$
   For $n\geq 1 $, simple estimations lead us to the following
    \[ \begin{split}
 \sum_{k=1}^{n+2}\frac{(n+3-k)^{\alpha-1}}{k}&=\left(\sum_{k=1}^{[\frac{n+2}{2}]}+\sum_{k=[\frac{n+2}{2}]+1}^{n+2}\right)\frac{(n+3-k)^{\alpha-1}}{k}\\
   & \asymp (n+1)^{\alpha-1}\sum_{k=1}^{[\frac{n+2}{2}]}\frac{1}{k} +\frac{1}{n+1} \sum_{k=[\frac{n+2}{2}]+1}^{n+2}(n+3-k)^{\alpha-1}\\
   & \asymp  (n+1)^{\alpha-1}\log(n+2)+(n+1)^{\alpha-1}\\
   & \asymp (n+1)^{\alpha-1}\log(n+2).
      \end{split} \]
      Therefore,
         \[ \begin{split}
         1&\gtrsim  \sn \mu_{n+2}\sum_{k=1}^{n+2}\frac{\Gamma(n+2-k+\alpha)}{\Gamma(\alpha)(n+2-k)!k}\\
         & \gtrsim  \sn (n+1)^{\alpha-1}\mu_{n}\log(n+2).
             \end{split} \]

 {\bf Case $p>1$.}  \ \ Let $q$ be the conjugate index of $p$, that is, $\frac{1}{p}+\frac{1}{q}=1$.  It is known that  $(B_{q})^{\ast}\cong B_{p}$(see \cite[Theorem 5.24]{H14}) under the paring
 $$\langle F,G\rangle=\int_{\dd}F'(z)\overline{G'(z)}dA(z), \ \ F\in B_{p},G\in B_{q}.$$
This means that $\cua$ is bounded from $X$ to $B_{p}$ if and only if
$$|\langle \cua(F),G\rangle|\lesssim ||F||_{X}||G||_{B_{q}}\ \mbox{for all}\ F\in X, G\in B_{q}.$$

Now, suppose $\cua$ is bounded from $X$ to $B_{p}$.  Take $g(z)=\sn\widehat{ g}(n)z^{n}\in B_{q}$  and the sequence of its Taylor coefficients is a
decreasing sequence of the non-negative real numbers. Let $f(z)=\log\frac{1}{1-z}=\sum_{n=1}\frac{z^{n}}{n}\in X$,  we have that
$$|\langle \cua(f),g\rangle|\lesssim ||f||_{X}||g||_{B_{q}}.$$
A simple calculation shows that
$$|\langle \cua(f),g\rangle|=\sum^{\infty}_{n=1}n\mu_{n}\left(\sum_{k=1}^{n}\frac{\Gamma(n-k+\alpha)}{\Gamma(\alpha)(n-k)!k}\right)\widehat{g}(n).$$
 This implies that
 $$|\langle \cua(f),g\rangle|=\sum^{\infty}_{n=1}n^{\frac{1}{q}}\mu_{n}\left(\sum_{k=1}^{n}\frac{\Gamma(n-k+\alpha)}{\Gamma(\alpha)(n-k)!k}\right)\widehat{g}(n)n^{\frac{q-1}{q}}<\infty.$$
 By Lemma \ref{lem2.1}, the sequence $\{\widehat{g}(n)n^{\frac{q-1}{q}}\}_{n=1}^{\infty}\in l^{q}$. The well known duality $(l^{q})^{\ast}=l^{p}$ yields that
 $$\left\{n^{\frac{1}{q}}\mu_{n}\left(\sum_{k=1}^{n}\frac{\Gamma(n-k+\alpha)}{\Gamma(\alpha)(n-k)!k}\right)\right\}^{\infty}_{n=1}\in l^{p}.$$
 Using the estimate $\sum_{k=1}^{n}\frac{\Gamma(n-k+\alpha)}{\Gamma(\alpha)(n-k)!k}\asymp (n+1)^{\alpha-1}\log(n+2)$  we deduce that
$$\sn (n+1)^{p\alpha-1}\mu_{n}^{p}\log^{p}(n+2)<\infty.$$
\end{proof of Theorem 1.1 step1}

\begin{proof of Theorem 1.1 step2}  
Let $\{f_{k}\}_{k=1}^{\infty}$   be a bounded sequence in  $X$  which converges to $0$ uniformly on every compact subset of  $\mathbb{D}$. Without loss of generality, we may assume   that $f_{k}(0)=0$ and $\sup_{k\geq 1}||f||_{X}\leq 1$. It suffices to prove that $\lim_{k\rightarrow \infty}||\cua(f_{k})||_{B_{p}}=0$
by using Lemma \ref{lem2.5}. As before, we divide the proof into $p = 1$ and $p >1$.

{\bf Case\textbf{ $p>1$.}} \ \ Assume  $\sum_{n=1}^{\infty}(n+1)^{p\alpha-1}\mu_{n}^{p}\log^{p}(n+1)<\infty$,  then
 \[ \begin{split}
\sum_{n=1}^{\infty}(n+1)^{p\alpha-1}\mu_{n}^{p}\log^{p}(n+1)&= \sum_{n=1}^{\infty}\left(\sum_{k=2^{n-1}}^{2^{n}-1}(k+1)^{p\alpha-1}\mu_{k}^{p}\log^{p}(k+1)\right)\\
& \gtrsim \sum_{n=1}^{\infty}2^{np\alpha}\mu_{2^{n}}^{p}\log^{p}(2^{n}+1)\\
& \gtrsim \sum_{n=1}^{\infty}2^{-n(p-1)}\left(\sum_{k=2^{n}}^{2^{n+1}-1}(k+1)^{\alpha-\frac{1}{p}}\mu_{k}\log(k+1)\right)^{p}.
  \end{split} \]
This shows that
$$ \sum_{n=1}^{\infty}2^{-n(p-1)}\left(\sum_{k=2^{n}}^{2^{n+1}-1}(k+1)^{\alpha-\frac{1}{p}}\mu_{k}\log(k+1)\right)^{p}<\infty.$$
By Lemma \ref{lem2.2} we have that
 \[ \begin{split}
& \ \ \ \  \int_{0}^{1}(1-r)^{p-2}\left(\sn (n+1)^{\alpha-\frac{1}{p}}\mu_{n}\log(n+1)r^{n}\right)^{p}dr \\
  &\asymp \sum_{n=0}^{\infty}2^{-n(p-1)}\left(\sum_{k=2^{n}}^{2^{n+1}-1}(k+1)^{\alpha-\frac{1}{p}}\mu_{k}\log(k+1)\right)^{p}<\infty.
   \end{split} \]
Therefore, for any  $\varepsilon>0$ there exists a  $0<r_{0}<1$ such that
\begin{equation}\label{1}
\int_{r_{0}}^{1}(1-r)^{p-2}\left(\sn (n+1)^{\alpha-\frac{1}{p}}\mu_{k}\log(k+1)r^{n}\right)^{p}dr<\varepsilon.
\end{equation}
It is clear that
\[ \begin{split}
||\cua(f_{k})||^{p}_{B_{p}}&=\left(\int_{|z|\leq r_{0}}+\int_{r_{0}<|z|<1}\right)|\cua(f_{k})'(z)|^{p}(1-|z|)^{p-2}dA(z)\\
& := J_{1,k}+J_{2,k}.
  \end{split} \]
By the integral representation of  $\cua$ we  get
  \begin{equation}\label{2}
  \Cu(f_{k})'(z)=\int_{0}^{1}\frac{tf'_{k}(tz)}{(1-tz)^{\alpha}}d\mu(t)+\int_{0}^{1}\frac{\alpha tf_{k}(tz)}{(1-tz)^{\alpha+1}}d\mu(t).\end{equation}
 Cauchy integral theorem implies that the sequence $\{f'_{k}\}_{k=1}^{\infty}$ is also converge to $0$ uniformly on every compact subset of  $\mathbb{D}$. Thus,  for  $|z|\leq r_{0}$ we have that
  \[ \begin{split}
|\cua(f_{k})'(z)|& \lesssim \int_{0}^{1}\frac{|f_{k}'(tz)|}{|1-tz|^{\alpha}}+ \frac{|f_{k}(tz)|}{|1-tz|^{\alpha+1}}d\mu(t)\\
& \lesssim \sup_{|w|<r_{0}}\left(|f_{k}(w)|+|f_{k}'(w)|\right)\int_{0}^{1}\frac{1}{(1-tr_{0})^{\alpha+1}}d\mu(t)\\
& \lesssim  \sup_{|w|<r_{0}}\left(|f_{k}(w)|+|f_{k}'(w)|\right).
    \end{split} \]
It follows that
$$J_{1,k}  \rightarrow 0, \ (k\rightarrow \infty).$$
Next, we estimate $J_{2,k}$.

Since $X\subset \mathcal {B}$,  we have
\begin{equation}\label{3}
|f_{k}(z)|\lesssim \log\frac{e}{1-|z|}\ \  \mbox{and}\  \ |f_{k}'(z)|\lesssim \frac{1}{1-|z|}\ \ \mbox{for all}\ \ k\geq 1, z\in \dd .\end{equation}
 By  (\ref{2}) and (\ref{3}), Minkowski inequity, Lemma \ref{lem2.4} we get
\[ \begin{split}
M_{p}(r,\cua(f_{k})') &= \left\{\int_{0}^{2\pi}\left|\int_{0}^{1}\frac{tf_{k}'(tre^{i\theta})}{(1-tre^{i\theta})^{\alpha}}+\frac{tf_{k}(tre^{i\theta})}{(1-tre^{i\theta})^{\alpha+1}}d\mu(t)\right|^{p}d\theta\right\}^{\frac{1}{p}}\\
&\lesssim  \left\{\int_{0}^{2\pi}\left(\int_{0}^{1}\frac{1}{(1-tr)|1-tre^{i\theta}|^{\alpha}}d\mu(t)\right)^{p}d\theta\right\}^{\frac{1}{p}}\\
& \ \ \ \ +\left\{\int_{0}^{2\pi}\left(\int_{0}^{1}\frac{\log\frac{e}{1-tr}}{|1-tre^{i\theta}|^{\alpha+1}}d\mu(t)\right)^{p}d\theta\right\}^{\frac{1}{p}}\\
&\lesssim \int_{0}^{1}\frac{1}{1-tr}\left(\int_{0}^{2\pi}\frac{d\theta}{|1-tre^{i\theta}|^{p\alpha}}\right)^{\frac{1}{p}}d\mu(t)\\
& \ \ \ \ +  \int_{0}^{1}\log\frac{e}{1-tr}\left(\int_{0}^{2\pi}\frac{d\theta}{|1-tre^{i\theta}|^{p(\alpha+1)}}\right)^{\frac{1}{p}}d\mu(t)\\
& \lesssim\int_{0}^{1}H(t,r)d\mu(t),
\end{split} \]
where
\[H(t,r)=\left\{
  \begin{array}{cc}
 \displaystyle{\frac{\log\frac{e}{1-tr}}{(1-tr)^{\alpha+1-\frac{1}{p}}},} & \text{if} \ \ \displaystyle{p>\frac{1}{\alpha}},\\
 \displaystyle{\frac{\log\frac{e}{1-tr}}{1-tr}, } &  \text{if} \ \ \displaystyle{p=\frac{1}{\alpha}}.
 \end{array}\right.\]
 Lemma \ref{lem2.3} yields that
 $$M_{p}(r,\cua(f_{k})') \lesssim\int_{0}^{1}H(t,r)d\mu(t) \asymp \sn (n+1)^{\alpha-\frac{1}{p}}\mu_{n}\log(n+1)r^{n}. $$
This together with  (\ref{1}) we have that
\[ \begin{split}
    J_{2,k}& =\int_{r_{0}<|z|<1}|\cua(f_{k})'(z)|^{p}(1-|z|)^{p-2}dA(z)\\
  &   \lesssim \int_{r_{0}}^{1}(1-r)^{p-2}M^{p}_{p}(r,\cua(f_{k})')dr\\
  & \lesssim   \int_{r_{0}}^{1}(1-r)^{p-2}\left(\sn (n+1)^{\alpha-\frac{1}{p}}\mu_{n}\log(n+1)r^{n}\right)^{p}dr\\
  & \lesssim \varepsilon.
    \end{split} \]
Consequently,
$$\lim_{k\rightarrow \infty}||\cua(f_{k})||_{B_{p}}=0.$$

{\bf Case\textbf{ $p=1$.}} \ \
When $p=1$, Lemma \ref{lem2.3} shows that the condition $\sn (n+1)^{\alpha-1}\mu_{n}\log(n+2)<\infty$ is equivalent to $\int_{0}^{1}\frac{\log\frac{e}{1-t}}{(1-t)^{\alpha}}d\mu(t)<\infty.$
Hence, for any $\varepsilon>0$ there exists a $0<t_{0}<1$ such that
\begin{equation}\label{4}
\int_{t_{0}}^{1}\frac{\log\frac{e}{1-t}}{(1-t)^{\alpha}}d\mu(t)<\varepsilon.
\end{equation}
By the integral representation of  $\cua$ we have
\begin{equation}\label{5}
\cua(f)''(z)=\int_{0}^{1}\left(\frac{t^{2}f''(tz)}{(1-tz)^{\alpha}}
+\frac{2\alpha t^{2}f'(tz)}{(1-tz)^{\alpha+1}}+
\frac{\alpha(\alpha+1)t^{2}f(tz)}{(1-tz)^{\alpha+2}}\right)d\mu(t).\end{equation}
For $0<r<1$, we have
\[ \begin{split}
M_{1}(r, \cua(f_{k})'')&\lesssim \sup_{|w|\leq t_{0}}\left(|f_{k}''(w)|+|f_{k}'(w)|+|f_{k}(w)|\right)\int_{0}^{t_{0}}\frac{d\mu(t)}{(1-t_{0}r)^{\alpha+2}}\\
& \ \ \ \ +\int_{0}^{2\pi}\int_{t_{0}}^{1}\frac{|f_{k}''(tz)|}{|1-tre^{i\theta}|^{\alpha}}+\frac{|f_{k}'(tz)|}{|1-tre^{i\theta}|^{\alpha+1}}
+\frac{|f_{k}(tz)|}{|1-tre^{i\theta}|^{\alpha+2}}d\mu(t)d\theta.
    \end{split} \]
    Since $\{f_{k}\}\subset X\subset \mathcal {B}$, we see that
   \begin{equation}\label{6}
   |f_{k}''(z)|\lesssim \frac{1}{(1-|z|)^{2}}\ \mbox{ for all }\  k \geq 1.\end{equation}
The assumption of $p$ means  that $\alpha \geq 1$. By Fubini theorem, (\ref{3}), (\ref{6}) and  Lemma \ref{lem2.4} we have
\[ \begin{split}
&\ \ \ \  \int_{0}^{2\pi}\int_{t_{0}}^{1}\frac{|f_{k}''(tre^{i\theta})|}{|1-tre^{i\theta}|^{\alpha}}+\frac{|f_{k}'(tre^{i\theta})|}{|1-tre^{i\theta}|^{\alpha+1}}
+\frac{|f_{k}(tre^{i\theta})|}{|1-tre^{i\theta}|^{\alpha+2}}d\mu(t)d\theta\\
&\lesssim \int_{t_{0}}^{1}\int_{0}^{2\pi}\left(\frac{1}{(1-tr)^{2}|1-tre^{i\theta}|^{\alpha}}+\frac{1}{(1-tr)|1-tre^{i\theta}|^{\alpha+1}}
+\frac{\log\frac{e}{1-tr}}{|1-tre^{i\theta}|^{\alpha+2}}\right)d\theta d\mu(t)\\
&\lesssim \int_{t_{0}}^{1} \frac{\log\frac{e}{1-tr}}{(1-tr)^{\alpha+1}}d\mu(t).
    \end{split} \]
Hence,
\[ \begin{split}
& \ \ \ \ \int_{t_{0}<|z|<1}|\cua(f_{k})''(z)|dA(z)\\
&= 2\int_{t_{0}}^{1}M_{1}(r, \cua(f_{k})'')rdr\\
& \lesssim \sup_{|w|\leq t_{0}}\left(|f_{k}''(w)|+|f_{k}'(w)|+|f_{k}(w)|\right) +\int_{t_{0}}^{1} \int_{t_{0}}^{1} \frac{\log\frac{e}{1-tr}}{(1-tr)^{\alpha+1}}d\mu(t)dr\\
& \lesssim \sup_{|w|\leq t_{0}}\left(|f_{k}''(w)|+|f_{k}'(w)|+|f_{k}(w)|\right) +\int_{t_{0}}^{1}\log\frac{e}{1-t}\int_{0}^{1}\frac{dr}{(1-tr)^{\alpha+1}}d\mu(t)\\
& \lesssim \sup_{|w|\leq t_{0}}\left(|f_{k}''(w)|+|f_{k}'(w)|+|f_{k}(w)|\right)+\int_{t_{0}}^{1}\frac{\log\frac{e}{1-t}}{(1-t)^{\alpha}}d\mu(t)\\
&\lesssim \sup_{|w|\leq t_{0}}\left(|f_{k}''(w)|+|f_{k}'(w)|+|f_{k}(w)|\right)+\varepsilon.
   \end{split} \]
The uniform convergence of $\{f_{k}\}$ on compact subset of $\dd$ implies that
$$\int_{|z|\leq t_{0}}|\cua(f_{k})''(z)|dA(z)\lesssim \sup_{|w|\leq t_{0}}\left(|f_{k}''(w)|+|f_{k}'(w)|+|f_{k}(w)|\right)\rightarrow 0, \ \ \mbox{as}\ k\rightarrow0. $$
Therefore, we deduce  that
$$\lim_{k\rightarrow \infty}||\cua(f_{k})||_{B_{1}}=0.$$
Thus, the operator $\cua$ is compact  from $X$  to $B_{1}$.
\end{proof of Theorem 1.1 step2}







\section*{Conflicts of Interest}
The authors declare that there is no conflict of interest.

\section*{Funding}
 The author was supported by  the Natural Science Foundation of Hunan Province (No. 2022JJ30369).


\section{Availability of data and materials}
Data sharing not applicable to this article as no datasets were generated or analysed during
the current study: the article describes entirely theoretical research.




 \end{document}